\documentclass[12pt]{article}

\usepackage[utf8]{inputenc}       
\usepackage[T1]{fontenc}          
\usepackage{lmodern}              
\usepackage{geometry}             
\usepackage{graphicx}             
\usepackage{amsmath, amssymb}     
\usepackage{hyperref}             
\usepackage[numbers]{natbib}      
\usepackage{setspace}             
\usepackage{booktabs}             
\usepackage{caption}              
\usepackage{subcaption}
\usepackage{pgfplots}
\usepackage{tikz}
\usepackage{algorithm}
\usepackage{algpseudocode}

\title{Concurrent Crossover for PDHG}
\author{Edward Rothberg \\
Gurobi Optimization, LLC \\
\texttt{rothberg@gurobi.com}}
\date{October 28, 2025}

\begin{document}

\maketitle
\begin{abstract}
  First-order methods based on the PDHG algorithm have recently
  emerged as a viable option for efficiently solving large-scale
  linear programming problems.  One highly desirable property of these
  methods is that they can make effective use of GPUs.  One
  undesirable property is that, as first-order methods, their
  convergence can be extremely slow.  This property forces one to
  decide how much accuracy is truly necessary when solving an LP
  problem.  This paper looks at whether a parallel, concurrent
  crossover scheme can help to obtain highly accurate solutions
  without sacrificing the benefits of these new approaches.
\end{abstract}

\section{Introduction}

A new option has emerged recently for solving large linear programming
problems.  A recent paper by a team at Google~\cite{pdlp22} proposed
an enhancement and specialization of an existing first-order method,
the Primal Dual Hybrid Gradient (PDHG) algorithm~\cite{pdhg11}, to the
problem of linear programming.  They call their approach {\em PDLP\/}.
The dominant computation in their algorithm, sparse matrix-vector
multiplication, is sufficiently simple and regular that it maps quite
well to massively parallel computer architectures, including GPUs.
Once GPU implementations of PDLP were built~\cite{cupdlp23} and
further enhancements were made~\cite{hprlp24,cupdlp+25}, this approach
was found to be competitive with the long-dominant approaches for
solving LPs, simplex~\cite{simplex51} and
interior-point~\cite{wright97}, when run on modern GPUs.

One issue with these new approaches is that, as first-order methods,
they can converge quite slowly to an optimal solution, leading to a
desire to terminate them as soon as practically possible.  As a
result, performance comparisons against simplex and interior-point
methods are typically not apples-to-apples, since PDHG methods are
typically stopped at a point where their solutions have errors that
are several orders of magnitude larger than those of the alternatives.

The disparity in the quality (and nature) of the solutions produced by
different LP methods isn't new.  Interior-point methods produce solutions
that are strictly interior to the LP feasible region, while simplex
produces {\em basic\/} solutions that are at a corner point.
Interior-point solutions also typically have substantially more error
(although not nearly as much as PDHG).  This disparity is typically
removed by employing a crossover method~\cite{crossover91}, which
efficiently converts an interior solution into a basic solution.
Unfortunately, in contrast to interior-point methods (and PDHG), this
crossover procedure has little ability to exploit parallelism.  As a
result, its relative cost compared with the interior-point solve is
large, and continues to grow as CPU core counts continue to increase.

The large errors that are typical in PDHG solutions can also be
resolved with a crossover approach, but that brings additional
challenges.  The time required for crossover depends heavily on the
quality of the starting point, so if PDHG iterations are stopped
early, crossover will likely take a long time.  Conversely, if PDHG is
allowed to proceed until it obtains a very accurate solution,
crossover will be faster but PDHG iterations may take a long time.

This paper looks at whether parallelism can be used to navigate this
tradeoff.  In particular, we consider an approach that launches
multiple concurrent crossover threads in parallel from multiple
intermediate PDHG solutions to accelerate the process of obtaining an
accurate, basic solution.  Measuring the impact presents a bit of
a challenge, but our results indicate that this approach yields an
average 25-50\% improvement in overall runtime.

\section{Background}

A linear programming problem can be stated in its primal and dual forms
as:
\begin{equation*}
\begin{minipage}{0.45\linewidth}
\begin{align*}
\min \quad & c^\top x \\
\text{s.t.} \quad & A x = b \\
& x \ge 0
\end{align*}
\end{minipage}
\hspace{0.2cm}
\begin{minipage}{0.45\linewidth}
\begin{align*}
\max \quad & b^\top y \\
\text{s.t.} \quad & A^\top y + z = c \\
                  & z \ge 0
\end{align*}
\end{minipage}
\end{equation*}
A candidate solution $(x,y,z)$ has an associated primal and dual residual
vector:
\begin{align*}
r_P &= b - Ax \\
r_D &= A^\top y + z - c
\end{align*}
A solution is optimal if (i) the primal and dual residual vectors are
zero, (ii) $x$ and $z$ are non-negative, and (iii) the objective gap is
zero.  This last condition can be expressed directly as
$c^\top x = b^\top y$, or somewhat less directly as $x^\top z = 0$
(referred to as {\em complementarity\/}).

Solutions produced by LP solvers always have some error, due to errors
in floating-point arithmetic but also due to a desire to limit
runtimes.  Different methods manage these errors differently.

The (primal) simplex method maintains primal feasibility and
complementarity (after an initial phase 1 to find a feasible basis),
and continues iterating until the maximum dual constraint violation is
below some threshold
($\|r_D\|_{\infty} \leq \varepsilon_{\text{abs}}$).  The traditional
choice is $\varepsilon_{\text{abs}} = 10^{-6}$.

PDHG maintains a current iterate that does not necessarily satisfy primal
feasibility, dual feasibility, or zero objective gap.  Iterations
proceed until all three achieve target tolerances (typically expressed
in relative terms):
\begin{align*}
  \| r_P \|_2 & \leq \varepsilon_{\text{rel}} (1 + \| b \|_2) \\
  \| r_D \|_2 & \leq \varepsilon_{\text{rel}} (1 + \| c \|_2) \\
  | c^\top x - b^\top y| & \leq \varepsilon_{\text{rel}} (1 + |c^\top x| + |b^\top y|)
\end{align*}
The most commonly used convergence tolerance when evaluating the
performance of these methods is $\varepsilon_{\text{rel}} = 10^{-4}$.
This is much weaker than the tighter, absolute tolerance typically used
for the simplex method.

Interior-point methods also maintain a current iterate that does not
necessarily satisfy feasibility or complementarity.  Interior-point
methods are {\em locally quadratically convergent\/}, meaning that
they converge quite quickly once the iterate is ``close enough'' to
an optimal solution.  Convergence theory puts strong bounds on the
number of iterations required to reach convergence, but these bounds
are based on relative tolerances.  Practice implementations typically
iterate well beyond the tolerances that convergence theory would
suggest.  For most models, these methods are able to achieve small absolute
violations ($\varepsilon_{\text{abs}} = 10^{-6}$), which gives
crossover a highly accurate starting solution.

Note that these different convergence criteria can have a profound
impact on solution quality.  To make this more concrete, consider the
well-known historical model {\em pilot87\/}, a modest-sized LP with
around 2,000 rows and 5,000 columns. PDHG with a $10^{-6}$ relative
tolerance on the 2-norm of the residual produces a solution with a
maximum absolute constraint violation of $4 \cdot 10^{-3}$.  The solution
found by the dual simplex method with a $10^{-6}$ absolute tolerance
on the infinity norm has a maximum violation of $5 \cdot 10^{-18}$.  One
can reasonably debate the question of whether a violation that small
is necessary, but the 15 order of magnitude difference will certainly
have an impact in some situations.

\subsection{Crossover}

Even setting aside issues of solution accuracy, not all optimal
$(x,y,z)$ are equally desirable.  A basic solution is the gold
standard in LP, for a number of reasons:
\begin{itemize}
\item Reoptimization: Simplex is quite effective at reoptimizing
  from a basic optimal solution after small changes to the model.
\item Sparsity: A basic solution puts most variables at bounds
  (which is particularly important for Mixed Integer Programming).
\item Conciseness: A solution can be captured to arbitrary precision
  by simply providing a list of the columns in the basis.
\end{itemize}
Given these advantages, it is not surprising that crossover has become
an integral ingredient for solving LP problems using interior-point
methods.  Even as the relative cost of crossover has grown
significantly, due to increasing use of parallelism for the
interior-point iterations, we have found it quite rare for users to turn
it off.

There is one notable exception, though, which is when crossover is
unable to reliably find an optimal basic solution, typically because
the interior-point solver struggles to produce an accurate starting
point.  While crossover implementations are built to be reliable in
the face of significant feasibility or complementarity violations,
there are limits to how much error they can accommodate.

\section{Crossover and PDHG}

On the surface, the question of when and how to apply crossover may
appear similar for interior-point and PDHG solvers.  There are some
important differences, though, driven by the difference in convergence
rates.  Interior-point solvers typically aim for violations that are as
small as possible, because smaller violations mean less time spent in
crossover.  This tradeoff was explored by Andersen and
Ye~\cite{multicross96}, where they considered the effectiveness of
starting crossover from earlier interior-point iterates.  While the
idea is intriguing, to our knowledge no modern implementations
incorporate it, likely because of the locally quadratically convergent
behavior of interior-point solvers.  If it typically only takes a few
iterations to evolve a solution that is close to an optimal solution
into a fully converged solution, the scope for improvement from
stopping early is small, and the risk of crossover time increasing
dramatically is large.

The linear convergence behavior of PDHG changes the dynamic
significantly.  Stopping PDHG iterations early and settling for a less
accurate crossover starting point can dramatically decrease the
convergence time, but it can significantly increase the crossover
time.  This difficult tradeoff motivates our concurrent crossover
approach.

\subsection{Concurrent Crossover Method}

To help identify a good PDHG iterate from which to launch crossover, our
proposed concurrent crossover PDHG (Algorithm 1) launches multiple crossover
threads from multiple iterates.
\begin{algorithm}
\caption{Concurrent Crossover PDHG}
\begin{algorithmic}[1]
    \State \textbf{Input:} Starting iterate $(x_0, y_0, z_0)$; $\varepsilon_{\text{cross}}$; $\varepsilon_{\text{rel}}$
    \State $\varepsilon_{\text{target}} \gets \varepsilon_{\text{cross}}$
    \While{$\text{maxresid}(x_k, y_k, z_k) > \varepsilon_{\text{rel}}$}
        \If{$\text{maxresid}(x_k, y_k, z_k) \leq \varepsilon_{\text{target}}$}
            \State Launch crossover thread from $(x_k, y_k, z_k)$ on the CPU
            \State $\varepsilon_{\text{target}} \gets 0.1 \, \varepsilon_{\text{target}}$
        \EndIf
        \State Perform PDHG iteration to produce $(x_{k+1}, y_{k+1}, z_{k+1})$
        \State $k \gets k + 1$
    \EndWhile
    \State Run crossover on main thread from converged iterate
\end{algorithmic}
\end{algorithm}
The method stops when the first crossover thread finds a basic,
optimal solution.

The approach requires a choice of target tolerance for the first
crossover thread ($\varepsilon_{\text{cross}}$), a choice of the
amount by which this violation must decrease to launch the next
crossover thread, and the tolerance at which the PDHG iterations cease
($\varepsilon_{\text{rel}}$).  A choice of
$\varepsilon_{\text{cross}}=10^{-2}$ and
$\varepsilon_{\text{rel}}=10^{-6}$, with a 10X reduction required
between threads, would launch up to four crossover threads, at $10^{-2}$,
$10^{-3}$, $10^{-4}$, and $10^{-5}$ maximum residuals.  It would also
launch crossover on the main thread once the target $10^{-6}$
tolerance is reached.  The first of these threads to find a basic
optimal solution returns its solution and terminates the others.

This approach is well suited to exploiting the characteristics of a
typical GPU system with dedicated GPU memory.  Nearly all PDHG
work is performed on the GPU and its memory system, leaving the CPU
nearly idle - one CPU thread launches kernels on the GPU and fetches
small amounts of data from GPU memory.  The only additional costs
imposed by the proposed concurrent crossover approach are: (1) the
cost of pulling the current iterate $(x_k,y_k,z_k)$ from GPU memory to
launch a new crossover thread, and (2) the contention for CPU
resources that comes with running multiple crossover threads on the
CPU.

Concurrent crossover imposes greater overheads when PDHG also runs on
the CPU, since threads launched for crossover are not available for
multi-threaded PDHG.  In addition, these threads compete for access to
the memory system with the main PDHG threads.  PDHG and crossover both
make heavy use of memory bandwidth, so this contention will slow down
all active threads on any modern system.  To manage the total thread
count, we simply reserve a few threads for crossover when we start
multi-threaded PDHG.  PDHG will saturate the memory system with a
modest number of threads, so removing a few typically does not
significantly hurt PDHG performance.  When crossover threads are
launched, they slow down the PDHG threads and any
previously-launched crossover threads.

\section{Results}

We will now evaluate the performance impact of this approach.

\subsection{Testing Environment}

All of our tests were performed using a pre-release version of Gurobi
version 13.0.  The PDHG implementation in this version closely follows
the methods described in cuPDLP+~\cite{cupdlp+25}.  Gurobi 13.0
includes both GPU and multi-threaded CPU implementations of PDHG.
Gurobi 13.0 also includes state-of-the-art implementations of the
primal simplex, dual simplex, interior-point, and crossover
algorithms, which have all been tested and tuned over many years and a
wide range of practical LP models.

When we compare against ``baseline'' PDHG, we use the Gurobi 13
default tolerance of $\varepsilon_{\text{rel}} = 10^{-6}$ to terminate
PDHG iterations and trigger crossover.  This value was empirically
chosen to provide a good balance between PDHG time and crossover time.

Our tests were performed on two platforms.  The CPU-only tests were
performed on a system with an AMD EPYC 7313P processor with 16 cores
and 128GB of memory.  This machine's memory system has a peak
bandwidth of 205~GB/s.  For concurrent crossover PDHG, we set
aside four cores for launching crossover threads.

GPU tests were performed on an Nvidia Grace Hopper GH200 system, which
consists of an Nvidia Grace CPU containing 72 cores and a Hopper H100
GPU.  The system has 512GB of shared memory, and the GPU has its own
96GB of HBM memory.  The GPU memory has a peak bandwidth of 4~TB/s.
Any direct comparisons between GPU and CPU performance in this paper
always come from this machine.

For reference, the Grace Hopper GPU performs PDHG iterations roughly
25 times faster than the AMD 7313P system for large models (and
roughly 11 times faster than the Grace Hopper CPU, but that comparison
is never relevant for the results presented here).  The Grace Hopper
CPU performs crossover iterations roughly 30\% faster than the AMD
CPU.

We used two sets of models for our testing.  The first is the
Mittelmann LP testset~\cite{hans25}, which is a set of 43 models that
are commonly used to benchmark various LP methods and implementations.
The models in this set are of moderate size and difficulty, with
solution times ranging from a few seconds to a few hundred seconds.
The majority solve in under 10 seconds.

Due to limitations with this set for our testing purposes, we also
collected a set of 90 models where PDHG shows signs of providing
performance advantages over simplex and interior-point.  These are
mostly drawn from proprietary customer models, so unfortunately we
are unable to make this set available.  We will refer to this set as the {\em
  PDHG friendly\/} set.  In general, the models in this set are much
larger and more difficult than the ones in the Mittelmann LP set.

We used a time limit of one hour for all runs.  We use shifted
geometric means~\cite{achterberg07} with a shift of 1 second when
comparing two methods over a set of models.

\subsection{Mittelmann LP testset}

Our first test looks at the performance benefit of our concurrent
crossover approach on the Mittelmann LP set.
The first line in
Table~\ref{tab:speedup_mittelmannlp}
shows shifted geometric means on the (AMD) CPU, comparing concurrent
crossover PDHG (using $\varepsilon_{\text{cross}} = 10^{-2}$ and
$\varepsilon_{\text{rel}} = 10^{-6}$)
against baseline PDHG.
Recall that the concurrent approach on the CPU
assigns 12 threads for PDHG, making 4 available for crossover.
The second line shows results on the GPU.
\begin{table}[htbp]
\centering
\begin{tabular}{c|c|c|c}
\textbf{} & \textbf{Performance Ratio} & \textbf{Wins} & \textbf{Losses}\\ \hline
CPU & 2.08 & 29 & 2 \\
GPU & 1.49 & 19 & 3 \\
\end{tabular}
\caption{Comparison of concurrent crossover PDHG against baseline PDHG, Mittelmann LP testset.}
\label{tab:speedup_mittelmannlp}
\end{table}
The table shows win and loss counts, where a method is considered to
have won (lost) if its runtime is at least 10\% smaller (larger) than
the baseline.

The results show a substantial advantage from the concurrent crossover
approach.  The CPU version produces more than a 2X mean speedup, with
better results on 29 models and worse on 2 (runtimes were within 10\%
on the other 12 models in the set).  The advantage on the GPU is
smaller (roughly 1.5X speedup), with fewer wins, but still a
substantial benefit.  It makes sense that the CPU improvement is
larger, because launching crossover earlier avoids PDHG iterations,
which are much more expensive on the CPU.

The few losses for concurrent crossover PDHG are somewhat surprising,
since the main thread performs the exact same computation as baseline
PDHG.  These losses are simply the result of the earlier, unproductive
crossover threads contending for resources with the winning main
thread.

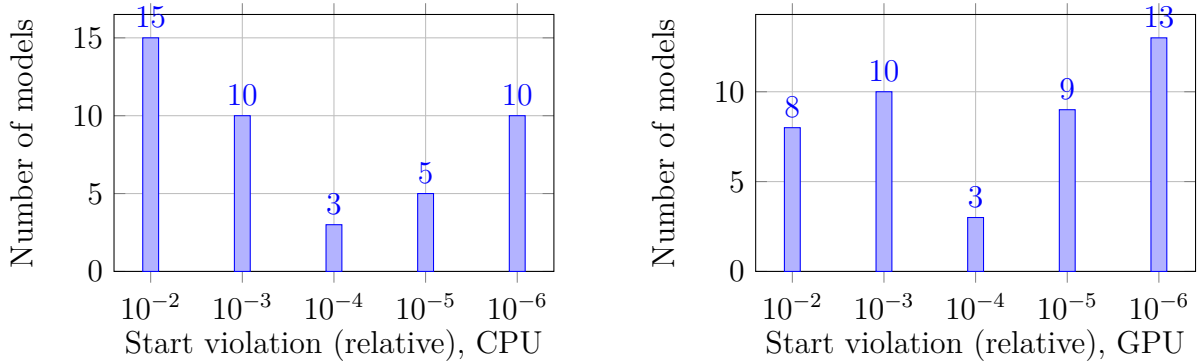
\begin{figure}[htbp]
\centering
\begin{subfigure}[t]{0.45\textwidth}
\centering
\begin{tikzpicture}
\begin{axis}[
    ybar,                                
    bar width=6pt,
    width=\linewidth,
    height=5cm,
    xlabel={Start violation (relative), CPU},
    xlabel style={yshift=-0.25em},
    ylabel={Number of models},
    symbolic x coords={$10^{-2}$,$10^{-3}$,$10^{-4}$,$10^{-5}$,$10^{-6}$}, 
    xtick=data,                          
    nodes near coords,                   
    ymin=0,                              
    grid=major,                          
]
\addplot coordinates {($10^{-2}$,15) ($10^{-3}$,10) ($10^{-4}$,3) ($10^{-5}$,5) ($10^{-6}$,10)};
\end{axis}
\end{tikzpicture}
\label{fig:win_frequency_mittelmannlp_cpu}
\end{subfigure}
\hspace{0.05\textwidth}
\begin{subfigure}[t]{0.45\textwidth}
\centering
\begin{tikzpicture}
\begin{axis}[
    ybar,                                
    bar width=6pt,
    width=\linewidth,
    height=5cm,
    xlabel={Start violation (relative), GPU},
    xlabel style={yshift=-0.25em},
    ylabel={Number of models},
    symbolic x coords={$10^{-2}$,$10^{-3}$,$10^{-4}$,$10^{-5}$,$10^{-6}$}, 
    xtick=data,                          
    nodes near coords,                   
    ymin=0,                              
    grid=major,                          
]
\addplot coordinates {($10^{-2}$,8) ($10^{-3}$,10) ($10^{-4}$,3) ($10^{-5}$,9) ($10^{-6}$,13)};
\end{axis}
\end{tikzpicture}
\end{subfigure}
\caption{Winning crossover thread for different starting violations (Mittelmann LP set).}
\label{fig:win_frequency_mittelmannlp}
\end{figure}

Figure~\ref{fig:win_frequency_mittelmannlp} shows how often each
crossover start point won the race to find an optimal basis.
You can see that for most of the 43 models in the testset, the
winning crossover thread on the CPU started from a PDHG solution with
either a $10^{-2}$ or $10^{-3}$ relative violation.
The winning thread on the GPU typically started with
a smaller violation, but $10^{-2}$ and $10^{-3}$ still win
quite often.

While these results are certainly encouraging,
Table~\ref{tab:speedup_mittelmannlp_concurrent} gives an indication
that they are missing some context.  This table shows the performance
benefit obtained over baseline PDHG by simply using the default LP
algorithm, which performs primal simplex, dual simplex, and
multi-threaded interior-point concurrently, terminating when
the first finishes.
\begin{table}[htbp]
\centering
\begin{tabular}{c|c|c|c|c}
\textbf{Platform} & \textbf{Method} & \textbf{Performance Ratio} & \textbf{Wins} & \textbf{Losses}\\ \hline
CPU & Concurrent crossover PDHG & 2.08 & 29 & 2 \\
    & Best of simplex and interior-point & 5.10 & 40 & 3 \\ \hline
GPU & Concurrent crossover PDHG & 1.49 & 19 & 3 \\
    & Best of simplex and interior-point & 2.75 & 21 & 6 \\
\end{tabular}
\caption{Comparison of default Gurobi against baseline PDHG, Mittelmann LP testset.}
\label{tab:speedup_mittelmannlp_concurrent}
\end{table}
Comparing these results against those for concurrent crossover
PDHG (included in this table as well) makes it clear that the
default Gurobi strategy is substantially faster.  In our
view, this significantly limits the value of evaluating our approach using
this set.  If all we have done is close a portion of the sizable
performance gap to the default approach, then we are not really
evaluating the practical impact of the approach.

\subsection{PDHG-friendly testset}

PDHG could certainly be incorporated into a concurrent scheme,
alongside the simplex and interior-point methods.  Concurrent
crossover would likely provide some benefit within that overall
scheme. However, that is an engineering question that we have chosen not
to tackle here.  Instead, we shift our attention to a more focused
question: would a concurrent crossover PDHG approach help for models
where the alternatives are known to be less effective?

As previously noted, we searched our LP test library for models where
PDHG was faster at finding a basic optimal solution than simplex or
interior-point methods, netting a set of 90 models.  The intent of
studying this (very biased) set is avoid situations where you can
reliably beat PDHG by simply not using PDHG.

Table~\ref{tab:speedup_vs_default_pdhgwin} shows results for the
Gurobi default LP solver (concurrent primal simplex, dual simplex, and
interior-point) compared against baseline PDHG on this testset,
on the CPU and the GPU.
\begin{table}[htbp]
\centering
\begin{tabular}{c|c|c|c}
\textbf{} & \textbf{Performance Ratio} & \textbf{Wins} & \textbf{Losses}\\ \hline
CPU & 1.32 & 42 & 33 \\
GPU & 0.33 & 4 & 82 \\
\end{tabular}
\caption{Comparison of default Gurobi against baseline PDHG, PDHG-friendly testset.}
\label{tab:speedup_vs_default_pdhgwin}
\end{table}
As you can see, baseline PDHG is 3X faster on these models on the GPU,
winning on all but a handful of models (the losses
can be explained by the fact that our search for our PDHG
friendly testset used a different Gurobi version).
Baseline PDHG is a bit behind on the CPU, but still quite close.

\begin{figure}[htbp]
\centering
\begin{subfigure}[t]{0.40\textwidth}
\begin{tikzpicture}
\begin{axis}[
    ybar,                                
    bar width=6pt,
    width=\linewidth,
    height=5cm,
    xlabel={Start violations (relative), CPU},
    xlabel style={yshift=-0.25em},
    ylabel={Number of models},
    symbolic x coords={$10^{-2}$,$10^{-3}$,$10^{-4}$,$10^{-5}$,$10^{-6}$}, 
    xtick=data,                          
    nodes near coords,                   
    ymin=0,                              
    grid=major,                          
]
\addplot coordinates {($10^{-2}$,12) ($10^{-3}$,20) ($10^{-4}$,17) ($10^{-5}$,31) ($10^{-6}$,10)};
\end{axis}
\end{tikzpicture}
\end{subfigure}
\hspace{0.05\textwidth}
\begin{subfigure}[t]{0.53\textwidth}
\begin{tikzpicture}
\begin{axis}[
    ybar,                                
    bar width=6pt,
    width=\linewidth,
    height=5cm,
    xlabel={Start violations (relative), GPU},
    xlabel style={yshift=-0.25em},
    ylabel={Number of models},
    symbolic x coords={$10^{-2}$,$10^{-3}$,$10^{-4}$,$10^{-5}$,$10^{-6}$,$10^{-7}$,$10^{-8}$}, 
    xtick=data,                          
    nodes near coords,                   
    ymin=0,                              
    grid=major,                          
]
\addplot coordinates {($10^{-2}$,0) ($10^{-3}$,6) ($10^{-4}$,15) ($10^{-5}$,23) ($10^{-6}$,23) ($10^{-7}$,12) ($10^{-8}$,11)};
\end{axis}
\end{tikzpicture}
\end{subfigure}
\caption{Winning crossover thread for different starting violations (PDHG-friendly testset).}
\label{fig:win_frequency_pdhgwin_gpu}
\end{figure}
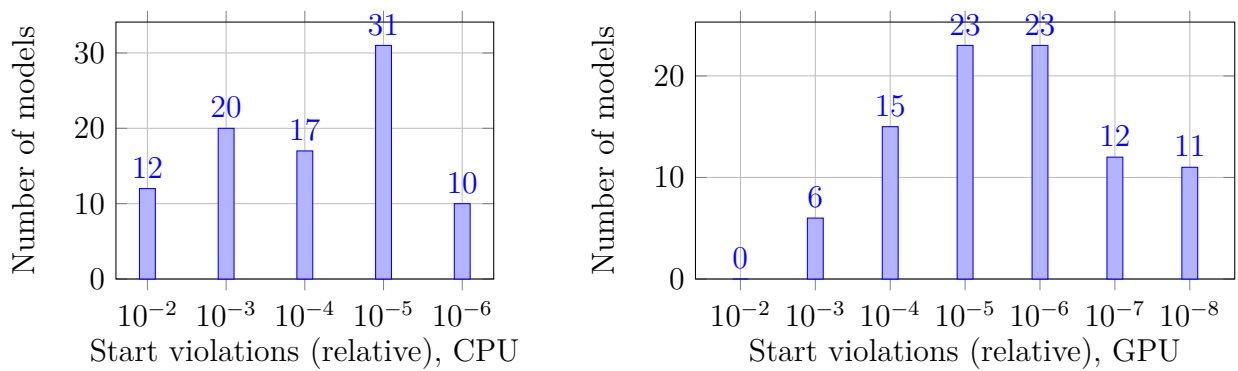

The histogram in Figure~\ref{fig:win_frequency_pdhgwin_gpu} shows how
often each crossover thread won.  We have increased the number of
concurrent crossover threads from 4 to 6 on the GPU (by choosing
$\varepsilon_{\text{rel}} = 10^{-8}$), to get a broader view of how
much accuracy is needed to get a good crossover result.  The data
shows that very loose tolerances rarely win on the GPU for these
harder models; the crossover threads we launch for these just get in
the way of the more productive threads.

Table~\ref{tab:speedup_pdhgwin} shows mean improvements on the CPU and
GPU from concurrent crossover PDHG.
\begin{table}[htbp]
\centering
\begin{tabular}{c|c|c|c}
\textbf{} & \textbf{Performance Ratio} & \textbf{Wins} & \textbf{Losses}\\ \hline
CPU ($\varepsilon_\text{cross} = 10^{-2}$) & 1.53 & 46 & 15 \\
GPU ($\varepsilon_\text{cross} = 10^{-2}$) & 1.20 & 32 & 22 \\
GPU ($\varepsilon_\text{cross} = 10^{-4}$) & 1.24 & 35 & 17 \\
\end{tabular}
\caption{Comparison of concurrent crossover PDHG against baseline PDHG, PDHG-friendly testset ($\varepsilon_\text{rel}=10^{-6}$ for the CPU and $10^{-8}$ for the GPU).}
\label{tab:speedup_pdhgwin}
\end{table}
We include two values of $\varepsilon_{\text{cross}}$ on the GPU,
$10^{-2}$ and $10^{-4}$, motivated by the fact that looser tolerances
rarely won.  The data shows that launching concurrent crossover
threads speeds up the overall time to a basic feasible solution
substantially - 1.53X on the CPU and 1.24X on the GPU.  While the
speedups are smaller than they were on the Mittelmann LP testset, the
results demonstrate that this approach does more than simply switch to
a faster alternative - it provides a boost even when PDHG is already the
clear winner.

\section{Discussion}

A natural question at this point is whether we are using the right
building blocks for finding a basic optimal solution from a PDHG
iterate.  One part of this is the question of whether a crossover
method that was designed and tuned for interior-point solvers is the
right tool if you expect the start point to have limited accuracy.
While interior-point solvers typically produce very accurate
solutions, they fail to fully converge fairly often, so lower accuracy
is something crossover already has to contend with.  There may be
better approaches for coping with this, but the underlying question has
already received quite a lot of attention.

Another possibility is that a completely new approach to crossover
could be developed that is better suited to parallel or GPU
architectures.  Liu and Lu~\cite{pdhgcrossover24} have considered a
specialized crossover approach inspired by PDHG, but this is early
work.

Another related question is whether PDHG could be modified to produce
an iterate that is easier to cross over from.  One option would be to
periodically change the objective function within the PDHG iterations
to coax the final iterate to a corner of the optimal face.  Such an
approach was explored by Ge, Wang, Xiong, and Ye~\cite{crossoverye25}
within the context of interior-point solvers, with some success.  For
PDHG, this would be similar to the {\em feasibility polishing\/} step
explored in~\cite{pdlp25}, where they periodically perform PDHG
iterations with a 0 objective instead of the original objective inside
the standard PDHG solve to focus attention on achieving primal
feasibility.  Our experiments with a similar approach that tries to
push the solution to a corner have shown limited success so far, but
we believe that there is scope here.

We suspect that the common narrative that PDHG is particularly
effective on massively parallel architectures is not strictly
accurate.  The bottleneck for PDHG is memory bandwidth, and while
massive bandwidth goes hand-in-hand with GPUs in current systems, our
understanding is that the two are not inextricably linked.  Similar
bandwidth could be provided on CPU-based machines if such machines
were commercially viable.  The Grace Hopper CPU, which is roughly 11X
slower than the GPU for PDHG iterations, can saturate its memory
system using just 10 of its 72 cores.  With a (much) higher-bandwidth
memory system, this same CPU could nearly keep pace with the GPU.
This may be an academic point, given the market reality that high
bandwidth is usually only available on GPU-based systems, but these
could be less tightly coupled in the future.

\section{Conclusion}

PDHG methods run on modern GPUs offer the potential to solve some
large, difficult LP problems faster than traditional methods.  They
force a new choice on users, though: do you want a fast solution or
an accurate solution?  While some application
domains can probably tolerate lower accuracy, we suspect that this
is not a welcome choice for most users.  Crossover offers a simple
way to avoid that choice, but it introduces another difficult
choice: when do you transition from PDHG iterations to crossover?
We have demonstrated that a concurrent crossover approach is
an effective way to exploit multiple CPU cores and typical GPU
machine architectures to substantially reduce the cost of
this crossover step.

\section*{Acknowledgements}

The author would like to thank Tobias Achterberg, Haihao Lu, Robert
Luce, and David Torres Sanchez for their helpful comments.

\bibliographystyle{plainnat}
\bibliography{references}

\end{document}